# Randomness as an equilibrium. Potential and probability density[1]

Marian Grendar, Jr.[*] and Marian Grendar[†]

[*]*Institute of Measurement Science, Slovak Academy of Sciences
Dúbravská cesta 9, Bratislava, 842 19, Slovakia. umergren@savba.sk*
[†]*marian.grendar@bb.telecom.sk*

**Abstract.** Randomness is viewed through an analogy between a physical quantity, density of gas, and a mathematical construct – probability density. Boltzmann's deduction of equilibrium distribution of ideal gas placed in an external potential field than provides a way of viewing probability density from a perspective of forces/potentials, hidden behind it.

## INTRODUCTION

That Probability and Mathematical Statistics can gain a new inspiration from Statistical Physics was already recognized and stressed by E. T. Jaynes. Jaynes's lifelong work was devoted in large extent to interplay of this three areas (see [3]).

In this peek we explore an analogy between gas density and probability density, mentioned at [1]. Boltzmann's deduction of equilibrium gas distribution in an external potential field than serves as a ground for viewing randomness from an equilibrium perspective.

Among the gains of such a weird effort one can find:

1) clear splitting of any random effect into its stochastic and deterministic parts/potentials,

2) physical interpretation of any probability distribution, through its relationship to potential. For instance, normal distribution is from the equilibrium perspective a distribution of linear harmonic oscillator (LHO). Ubiquity of LHO in the world of Physics and ubiquity of normal distribution in the world of Probability and Statistics are thus related,

3) a support for search for causes/forces in study of a random effect. Revealing forces behind a random effect thus can help to understand the effect, and consequently help to predict (or in the case of luck also govern/regulate) the future outcomes of the effect,

4) an interpretation of important yet (almost) meaningless statistical quantities, like Shannon's entropy or Fisher's information. Also, MaxEnt has a nice interpretation in the presented context.

## RANDOMNESS AS AN EQUILIBRIUM

### The fundamental equivalence

Probability mass/density function is a mathematical construct. Statistical Physics offers a physical model of the probability distribution – an ideal gas, placed in an external potential field. Probability density function is materialized by the density function of the gas, since probability of finding a molecule of gas at an elementary volume is identical with the number of molecules of the gas at the volume.

The equivalence of the gas and probability densities lays in the ground of the presented view of randomness. Its equilibrium attribute then comes from Boltzmann's deduction of equilibrium distribution of the ideal gas placed in an external *potential field*.

### Boltzmann's formula

Boltzmann's deduction of equilibrium distribution of ideal gas placed in an external potential field, in a textbook exposition (see for instance [4]), can be for our purposes given as follows:

---

Molecules of the gas, closed in an infinitesimally small cube with sides $dx, dy, dz$ are subjected to two forces: force $F^1$, represented also by a potential function $U(x, y, z)$, and force $F^2$, induced by spatial difference of density $n(x, y, z)$ of the gas. The two forces are given (restricting to one-dimensional case only) by

$$dF^1 = -\frac{dU(x)}{dx} n(x) \, dx$$

and

$$dF^2 = -\frac{dn(x)}{dx} \, dx$$

In the *equilibrium*, the two forces should compensate each other, so

$$dF^1 + dF^2 = 0$$

hence

$$\frac{1}{n(x)} \frac{dn(x)}{dx} = -\frac{dU(x)}{dx}$$

The equilibrium condition is a differential equation, solved by

$$n(x) = ke^{-U(x)}$$

– the famous Bolzmann's formula. The formula expresses distribution of the ideal gas placed in a potential field, *in the equilibrium*.

## Basics

Recalling the above mentioned fundamental equivalence of the gas and probability densities, Boltzmann's deduction can be transferred into pure probability, and put in the ground of its equilibrium view.

Principles of the physicalization of randomness, admittedly metaphysical, can be formulated as follows:

$P1$ : Randomness has its own, *intrinsic* field.

$P2$ : Any randomness takes place also in an external field.

Fundamental terms of vocabulary of the equilibrium view of randomness are built up in analogy with the above mentioned Boltzmann's deduction.

**Definition 1** *Let X be a random variable, with pdf/pmf $f_X(x)$ defined over a support S. Intrinsic intensity $E^s$ of its potential field is defined as*

$$E^s(x) = -\frac{1}{f_X(x)} \frac{df_X(x)}{dx} \quad (1)$$

$E^s$ will be also called *stochastic component of randomness, or stochastic intensity*.

*Intensity of external field is*

$$E^c(x) = -\frac{dU(x)}{dx} \quad (2)$$

*where $U(x)$ is potential function of the external field. $E^c$ will be also called causal component of randomness, or causal intensity.*

**Definition 2** *If intensities of the internal and external fields compensate each other, the random variable attains its equilibrium distribution in potential $U(x)$*

$$f_X(x) = ke^{-U(x)} \quad (3)$$

*where $k = 1/\Omega$ is normalizing constant, $\Omega = \sum e^{-U(x)}$ is the statistical sum.*

**Definition 3** *Normalized potential $\tilde{U}(x)$ of the random variable X is*

$$\tilde{U}(x) = U(x) - \ln(k) \quad (4)$$

**Note 1** *According to (3) and (4) there is a unique relationship between normalized potential and equilibrium distribution:*

$$\tilde{U}(x) = -\ln(f_X(x)) \quad (5)$$

## Causal Intensity and Equilibrium Probability

According to defining formulas (2), (3), (4) and (5) equilibrium probability distribution generated by causal intensity $E^c(x)$, which is in the equilibrium in absolute value identical with the stochastic intensity $E^s(x)$, is

$$f_X(x) = e^{\int E^c(x)dx + c}$$

where $c$ is the integration/normalization constant.

In this section we will discuss the simplest form causal intensities and their corresponding equilibrium distributions, together with their normalized potentials.

1) *Zero causal intensity* $E^c(x) = 0$ generates *uniform distribution*. Equivalently: uniform distribution reveals absence of causal component of randomness. Uniform distribution is intrinsic distribution of randomness, left alone. Normalized potential of the uniform $N$-element discrete distribution is $\tilde{U}(x) = \ln(N)$.

2) *Constant causal intensity* $E^c(x) = -a$ generates *exponential distribution* $Exp(a)$. Equivalently: exponential distribution reveals presence of constant causal intensity in a studied effect. Normalized potential of the exponential distribution is $\tilde{U}(x) = ax - ln(a)$.

3) *Linear causal intensity* $E^c(x) = -bx$ generates *normal distribution* $n(0, \sqrt{\frac{1}{b}})$. Equivalently: normal distribution reveals a linear causal intensity. Normalized potential of $n(\mu, \sigma)$ is a *linear harmonic oscillator* potential $\tilde{U}(x) = \frac{(x-\mu)^2}{2\sigma^2} + \ln(\sqrt{2\pi\sigma^2})$.

4) *Superposition of linear and constant intensities* $E^c(x) = -a - bx$ generates distribution $f_X(x) = ke^{-ax - b\frac{x^2}{2}}$ where $k$ is a normalizing constant.

5) *Causal intensity of the form* $E^c(x) = -\frac{\Gamma'(x+1)}{\Gamma(x+1)} + \ln(\lambda)$ generates Poisson distribution $Poi(\lambda)$. Equivalently, Poisson distribution reveals presence of the above force. Normalized potential of Poisson distribution is $\tilde{U}(x) = \lambda - \ln(\lambda)x + \ln(\Gamma(x+1))$.

6) *Causal intensity* $E^c(x) = -\frac{a}{x} - b$ generates gamma distribution $\Gamma(1-a, \frac{1}{b})$. Equivalently: gamma distribution reveals presence of a superposition of constant and reciprocal causal intensities. Normalized potential of $\Gamma(\alpha, \beta)$ distribution is $\tilde{U}(x) = (1-\alpha)\ln(x) + \frac{x}{\beta} + \ln(\Gamma(\alpha)\beta^\alpha)$.

7) Pearson's system of distributions, an appreciated tool for various densities approximation (see [5]), is defined by differential equation

$$\frac{1}{f(x)} \frac{\partial f(x)}{\partial x} = \frac{x-a}{b_0 + b_1 x + b_2 x^2}$$

As it can be seen from (2), (4) and (5), the right-hand side of the differential equation expresses causal intensity $E^c(x)$ which generates the Pearson's system.

## MaxEnt: equilibrium of causal and stochastic intensities

First, we recall a definition of MaxEnt task (see [1]) simplified, without loss of generality, to one-potential-function case:

**Definition 4** *ME task with $u(x)$. Given a random sample and a known potential function $u(x)$, the maximum entropy task is to find the most entropic distribution **p** consistent with u-moment consistency condition.*

There, we made use of $u$-moment consistency condition:

**Definition 5** *Let $\mu(u)$, $m(u)$ are u-moment and sample u-moment, respectively. Then requirement of their equality*

$$\mu(u) = m(u)$$

*is called u-moment consistency condition.*

Also, the following well-known result (see for instance [1]) should be recalled: The most entropic distribution **p** satisfying the $u$-moment consistency condition is the simple exponential form pmf/pdf $f_X(x|\lambda) = k(\lambda)e^{-\lambda u(x)}$, where $\lambda$ is such that $u$-potential moment consistency condition is satisfied.

**Theorem 1** *MaxEnt task with $u(x)$-potential is solved by* equilibrium distribution *which has normalized potential $\tilde{U}(x) = \lambda u(x) - \ln(k)$.*

Thus, the presented view of randomness offers an interpretation of Shannon's entropy maximization as a search for such a probability distribution, which balances intrinsic field (stochastic component of randomness) with given external potential field $u(x)$ (causal component of randomness), multiplied by a multiplicative constant $\lambda$. The presence of $\lambda$ can appear for the first glance disturbing, but – as we claim at [2] – it is just the distinctive feature of Shannon's entropy maximization, and as such can give an intuitive answer to eternal 'Why MaxEnt?' question (for the argument, please see [2]).

It is also worth noting, that Shannon's entropy itself can be interpreted in the proposed context as *mean normalized potential*. Fisher's information number for the simple exponential form distribution has also very nice meaning – it is variance of potential.

## ACKNOWLEDGEMENTS

It is a pleasure to thank Aleš Gottvald, George Judge and Viktor Witkovský for valuable discussions and/or comments on earlier version of this work. The thanks extend to participants of MaxEnt workshop, especially to Ariel Caticha, Peter Cheeseman, Robert Fry, Samuel Kotz, Carlos Rodríguez and Alberto Solana for their interest and/or comments to presented work. Carlos Rodríguez pointed us to an interesting work by Kenneth Hanson and Gregory Cunningham (see [6]). The work was in part supported by the grant VEGA 1/7295/20 from the Scientific Grant Agency of the Slovak Republic.

## REFERENCES

1. M. Grendár, Jr., and M. Grendár, "MiniMax Entropy and Maximum Likelihood: Complementarity of Tasks, Identity of Solutions," in *Bayesian Inference and Maximum Entropy Methods in Science and Engineering*, edited by A. Mohammad-Djafari, pp. 49–61, AIP Press, New York, 2001. Available on-line at http://xxx.lanl.gov/abs/math.PR/0009129.
2. M. Grendár and M. Grendár, Jr., "Why Maximum Entropy? A Non-axiomatic Approach," in *Bayesian Inference and Maximum Entropy Methods in Science and Engineering*, edited by R. Fry, AIP Press, New York; this volume.
3. E. T. Jaynes, *Papers on Probability, Statistics and Statistical Physics*, edited by R.D. Rosenkrantz, D. Reidel Publishing Co., Dordrecht, 1983. Also, visit bayes.wustl.edu, web-page maintained by Larry Bretthorst.
4. A. N. Matveyev, *Molecular Physics*, Vysshaya Shkola, Moscow, 1987 (in Russian).
5. A. Stuart and J. K. Ord, *Kendall's advanced theory of Statistics, Vol. 1, Distribution theory*, Edward Arnold, London, 1994.

6. K. M. Hanson and G. S. Cunningham, "The Hard Truth," in *Maximum Entropy and Bayesian Methods*, edited by J. Skilling and S. S. Sibisi, pp. 157-164, Kluwer Academic, Dordrecht, 1996.